\documentclass[12pt,a4paper]{article}
\usepackage{amsmath,amsfonts,amssymb,amsthm}
\usepackage{mathrsfs}
\usepackage{ifthen}

\newtheorem{Theorem}{Theorem}
\newtheorem{Cor}{Corollary}
\newtheorem{Def}{Definition}
\newtheorem{Lemma}{Lemma}
\newtheorem{Prop}{Proposition}

\makeatletter
\newcommand\txt[2]{{#1{\ifthenelse{\equal {#2}{}}{}{ #2}}. }}
\def\Zam#1.#2\par{\txt{\bf Замечание}{#1}\ #2\par%
               \@ifnextchar {\Zam}{}{\@ifnextchar\begin{}{\smallskip}}}
\def\Example#1.{\txt{\bf Пример}{#1}}
\makeatother

\newcommand\emp\varnothing
\newcommand\eps\varepsilon
\newcommand\ov\widetilde
\let\ds\displaystyle

\let\dd\partial

\def\dout{\dd^{\rm out}}
\def\dist{\mathop{\operatorfont dist}}
\def\^#1{^{\overline{#1}}}

\begin{document}

\title{On graphs with a large chromatic number containing no small odd cycles}

\author{S.L. Berlov\thanks{The work was supported by RFBR grant No.~10-01-00096-A.}, Ilya I. Bogdanov\footnotemark[1]}

\maketitle

\begin{abstract}
  In this paper, we present the lower bounds for the number of vertices in a graph with a large chromatic number containing no small odd cycles.
\end{abstract}

\section{Introduction}

P.~Erd\H os~\cite{erdos} showed that for every integer $n>1$ and $p>2$, there exists a graph of girth~$g$ and chromatic number greater than~$n$ which contains not more than $n^{2g+1}$ vertices. Later, he conjectured~\cite{erdos3} that for every positive integer~$s$ there exists a constant $c_s$ such that for every graph~$G$ having $N$ vertices and containing no odd cycles of length less than $c_sN^{1/s}$, its chromatic number does not exceed $s+1$.

This conjecture was proved by Kierstead, Szemer\'edi, and Trotter~\cite{kierstead}; in fact, they have proved a more general result. In our case, their result states that the chromatic number of any graph on $N$ vertices containing no odd cycles of length at most $4sN^{1/s}+1$ does not exceed~$s+1$.

Basing on these results, we introduce the following notation.

\begin{Def}
  Assume that $n,k>1$ are two integers. Denote by $f(n,k)$ the maximal integer $f$ satisfying the following property: If a graph $G=(V,E)$ contains no odd cycles of length at most $2k-1$, and $|V|\leq f$, then there exists a proper coloring of its vertices in $n$ colors.
\end{Def}

Notice that a graph contains no odd cycles of length at most $2k-1$ if and only if it contains no simple odd cycles of the same lengths.

The results mentioned above imply that  $f(n,k)<n^{4k+1}$ and
$f({s+1},[2sN^{1/s}]+1)\geq N$. One can obtain that the latter inequality is equivalent to the bound
\begin{equation}
  f(n,k)\geq \left(\frac k{2(n-1)}\right)^{n-1}-1.
  \label{lower-kst}
\end{equation}

A different upper bound for $f(n,k)$ can be obtained from the following graph constructed by Schrijver~\cite{schrijver}. Let $m,d$ be some positive integers. Set $X=\{1,2,\dots,2m+d\}$, $V=\{{x\subset X}:\; |x|=m,\; 1<|i-j|<{2m+d-1}
\mbox{ \ for all pairs of distinct\ }i,j\in x\}$, $E=\{{(x,y)\in V^2}:\; {x\cap y=\emp}\}$. The Schrijver graph $(V,E)$
is $(d+2)$-chromatic, whilst it does not contain odd cycles of length less than $\frac{2m+d}d$. Next, we have $|V|=\frac{2m+d}{m+d}\binom{m+d}{d}$; now it is easy to obtain that
\begin{equation}
  f(n,k)<\frac{(n-1)(2k-1)+2}{(n-1)k+1}\binom{(n-1)k+1}{n-1}.
  \label{upper-schr}
\end{equation}
When we fix the value of $n$, the bounds~\eqref{lower-kst} and~\eqref{upper-schr} become the polynomials in $k$ of the same degree; hence, in some sense they are close to each other. On the contrary, when we fix the value of $k$ and consider the values $n>k/(2e)+1$, we see that the right-hand part of~\eqref{lower-kst}
decreases (as a function in $n$). Hence for larger values of $n$ this estimate does not provide any additional information.

On the other hand, for $k=2$ the asymptotics of $f(n,2)$ is tightly connected with the asymptotics of Ramsey numbers $R_{n,3}$. In the papers of Ajtai, Koml\'os, and E. Szemer\'edi~\cite{aks} and Kim~\cite{kim} it is shown that $c_1\frac{n^2}{\log n}\leq
R(n,3)\leq c_2\frac{n^2}{\log n}$ for some absolute constants $c_1,c_2$. One can check that these results imply the bounds
$$
  c_3 n^2\log n\leq f(n,2)\leq c_4 n^2 \log n
$$
for some absolute constants $c_3,c_4$.

In the present paper, we find nontrivial lower bounds for $f(n,k)$ for all values of $n\geq 2$ and $k\geq 2$. In Section~2, we make some combinatorial considerations leading to the recurrent bounds for $f(n,k)$. In Section~3, we obtain explicit bounds following from those results. In particular, we show (see Theorem~\ref{estimate}) that
$$
  f(n,k)\geq \frac{(n+k)(n+k+1)\cdots(n+2k-1)}{2^{k-1}k^k}
$$
for all $n\geq2$ and $k\geq 2$.

\section{Recurrent bounds}

Firstly, we introduce some notation.

Let $G=(V,E)$ be an (unoriented) graph. We denote the {\em distance} between the vertices $u,v\in V$ by $\dist_G(u,v)$.

Consider a vertex $v\in V$, and let $r$ be a nonnegative integer. We denote by $U_r(v,G)=\{u\in V\mid \dist_G(u,v)\le r\}$ the {\em ball} of radius $r$ with the center at~$v$, and by $S_r(v,G)=\{u\in V\mid \dist_G(u,v)= r\}$ the {\em sphere} with the same radius and center. In particular, $S_0(v,G)=U_0(v,G)=\{v\}$. Denote also by $\dout_G
V_1=\{u\in V\setminus V_1\mid
\exists v\in V_1: (u,v)\in E\}$ the {\em outer boundary} of a set $V_1\subseteq V$. In particular, $S_r(v,G)=\dout_G U_{r-1}(v,G)$.

For a set $V_1\subseteq V$, we denote by $G(V_1)$ the induced subgraph on the set of vertices~$V_1$.

\smallskip
Let us fix some integers $n$ and $k$ which are greater than~1. We need the following easy proposition.

\begin{Prop}
  Graph $G$ does not contain odd cycles of length not exceeding $2k-1$ if and only if for each vertex $v\in V$ and each positive integer $r<k$, the subgraph $G(S_r(v,G))$ contains no edges.
  \label{sph-empty}
\end{Prop}

\proof Assume that the subgraph $G(S_r(v,G))$ contains an edge $(u_1,u_2)$. Supplementing this edge by shortest paths from $v$ to $u_1$ and $u_2$, we obtain a cycle of length~$2r+1\leq 2k-1$.

Conversely, assume that $G$ contains a cycle of length $\leq 2k-1$.
Consider such a cycle $C$ of the minimal length $2r+1$ (then $r<k$). Choose any its vertex~$v$, and let $u_1,u_2$ be two vertices of~$C$ such that $\dist_C(v,u_1)=\dist_C(v,u_2)=r$. In fact, we have
$\dist_G(u_1,v)=\dist_G(u_2,v)=r$. Actually, assume that $\dist_G(v,u_1)<r$, and choose a path~$P$ of the minimal length connecting $v$ and $u_1$. Then one can supplement it by one of the two subpaths of~$C$ connecting $u_1$ and $v$ to obtain an odd cycle $C'$. The length of $C'$ is smaller than $r+(r+1)=2r+1$, that contradicts the choice of~$C$.

Thus, $u_1,u_2\in S_r(v,G)$, and the graph $G(S_r(v,G))$ contains an edge.
\qed

Now let us fix an arbitrary graph $G=(V,E)$ with a minimal number of vertices such that it contains no odd cycles of length not exceeding  $2k-1$, and $\chi(G)>n$ (hence $|V|=f(n,k)+1$). By the minimality condition, the graph~$G$ is connected. Moreover, for every $v\in V$ and $0\leq r\leq k$, the sphere~$S_r(v,G)$ is nonempty. Otherwise we would have $G=\cup_{i=0}^{r-1} S_i(v,G)$, where all the graphs $G(S_i(v,G))$ contain no edges  by Proposition~\ref{sph-empty}. Therefore, it is possible to color this graph properly in two colors: the vertices of the sets $S_i(v,G)$ with even~$i$ in color~1, while those for odd~$i$ --- in color~2 (the vertex~$v$ should be colored in color~1).

\smallskip
Let us introduce the number $d=\max_{v\in V} |U_{k-1}(v,G)|$.

\begin{Lemma}
\label{big okr}
  For every vertex $v\in V$, we have $|U_{k-1}(v,G)|\geq n(k-1)+1$. In particular, $d\geq n(k-1)+1$.
\end{Lemma}

\proof
Notice that $U_{k-1}(v,G)=\bigcup_{r=0}^{k-1}S_r(v,G)$. Assume that $|S_r(v,G)|\geq n$ for every $r=1,\dots,k-1$; then
$$
  |U_{k-1}(v,G)|=\sum_{r=0}^{k-1}|S_r(v,G)|\geq 1+(k-1)\cdot n,
$$
as desired.

Assume now that $|S_r(v,G)|< n$ for some $1\leq r\leq k-1$. Consider a subgraph $G'=G\bigl(V\setminus U_{r-1}(v,G)\bigr)$. From the minimality condition, it can be colored properly in $n$ colors. Consider an arbitrary such proper coloring; then the vertices of $S_r(v,G)$ are colored in at most $n-1$ colors, so there exists a color (say, color~1) different from them. Let us now color the vertices of~$S_{r-1}(v,G)$ in color~1, and then color all the remaining vertices of the sets~$S_i(v,G)$ ($i<r-1$) alternately: we use colors~1 and~2 (here 2 is any remaining color) for even and odd values of~$i-(r-1)$, respectively. It follows from Proposition~\ref{sph-empty} that this coloring is proper. This contradicts the choice of~$G$.
\qed

\begin{Lemma}
  \label{Sergey}
  $|V|\geq f(n-1,k)+d+1$.
\end{Lemma}

\proof
Choose a vertex~$v$ such that $d=|U_{k-1}(v,G)|$. Assume that $|V\setminus
U_{k-1}(v,G)|\leq f({n-1},k)$; then one can color properly vertices of the set $V\setminus U_{k-1}(v,G)$ in $n-1$ colors. Now we can color the vertices of the set $U_{k-1}(v,G)$ in colors~1 and~$n$ (where $n$ is a new color, and 1 is any of the old colors) in the following way: we color all the vertices of~$S_r(v,G)$ in color~1 or~$n$, if $r-(k-1)$ is odd or even, respectively. By Proposition~\ref{sph-empty}, we obtain a proper coloring of~$G$ in $n$ colors which is impossible.

Thus, our assumption is wrong, so $|V\setminus
U_{k-1}(v,G)|\geq f(n-1,k)+1$, and
$$
  |V|\geq f(n-1,k)+1+|U_{k-1}(v,G)|=f(n-1,k)+d+1.
  \eqno\qed
$$

\begin{Lemma}
  \label{Ilya}
  $\ds |V|\geq \frac{d^{1/(k-1)}}{d^{1/(k-1)}-1}\bigl(f(n-2,k)+1\bigr)$.
\end{Lemma}

\proof
We will construct inductively a sequence of partitions of~$V$ into nonintersecting parts,
$$
  V=U_1\sqcup U_2\sqcup\dots\sqcup U_s\sqcup N_s\sqcup V_s,
$$
such that the following conditions are satisfied:

(i) for all $i=1,\dots,s$ we have $\dout_G U_i\subseteq N_s$; moreover, $\dout_G V_s\subseteq N_s$;

(ii) for every $i=1,2,\dots,s$ the graph $G(U_i)$ is bipartite (in fact, $U_i$ is a ball with radius not exceeding~$k-1$ in a certain subgraph of~$G$);

(iii) $(d^{1/(k-1)}-1)(|U_1|+\dots+|U_s|)\geq |N_s|$.

For the base case $s=0$, we may set $V_0=V$, $N_0=\emp$ (there are no sets~$U_i$ in this case).

For the induction step, suppose that the partition $V=U_1\sqcup U_2\sqcup\dots\sqcup U_{s-1}\sqcup N_{s-1}\sqcup V_{s-1}$ has been constructed, and assume that the set~$V_{s-1}$ is nonempty. Consider the graph $G_{s-1}=G(V_{s-1})$ and choose an arbitrary vertex $v\in V_{s-1}$. Now consider the sets
$$
  U_0(v,G_{s-1})=\{v\},\quad  U_1(v,G_{s-1}), \quad \dots, \quad U_{k-1}(v,G_{s-1}).
$$
One of the ratios
$$
  \frac{|U_1(v,G_{s-1})|}{|U_0(v,G_{s-1})|}, \quad \frac{|U_2(v,G_{s-1})|}{|U_1(v,G_{s-1})|}, \quad
  \dots, \quad \frac{|U_{k-1}(v,G_{s-1})|}{|U_{k-2}(v,G_{s-1})|}
$$
does not exceed~$d^{1/(k-1)}$, since the product of these ratios is
$$
  |U_{k-1}(v,G_{s-1})|\leq |U_{k-1}(v,G)|\leq d.
$$
So, let us choose $1\leq m\leq k-1$ such that
$$
  \frac{|U_m(v,G_{s-1})|}{|U_{m-1}(v,G_{s-1})|}\leq d^{1/(k-1)}.
$$

Now we set
$$
  U_s=U_{m-1}(v,G_{s-1}), \quad N_s=N_{s-1}\cup S_m(v,G_{s-1}),
  \quad
  V_s=V_{s-1}\setminus U_m(v,G_{s-1}).
$$
Since the condition~(i) was satisfied on the previous step, we have
$$
  \dout_G V_s\subseteq \dout_G V_{s-1}\cup S_m(v,G_{s-1})\subseteq N_s
$$
and
$$
  \dout_G U_s\subseteq \dout_G V_{s-1}\cup S_m(v,G_{s-1})\subseteq N_s,
$$
so this condition also holds now. The condition~(ii) is satisfied by Proposition~\ref{sph-empty}. Finally, the choice of~$m$ and the condition~(iii) for the previous step imply that
\begin{gather*}
  d^{1/(k-1)} |U_s|=d^{1/(k-1)}|U_{m-1}(v,G_{s-1})|\geq |U_m(v,G_{s-1})|,\\
  (d^{1/(k-1)}-1)(|U_1|+\dots+|U_{s-1}|)\geq |N_{s-1}|
\end{gather*}
and hence
$$
  (d^{1/(k-1)}-1)(|U_1|+\dots+|U_{s}|)\geq
  |N_{s-1}|+|U_m(v,G_{s-1})|-|U_{m-1}(v,G_{s-1})|
  =|N_s|.
$$
Thus, the condition~(iii) also holds on this step.

Continuing the construction in this manner, we will eventually come to the partition with $V_s=\emp$ since the value of $|V_s|$ strictly decreases. As the result, we obtain the partition $V=U_1\sqcup U_2\sqcup\dots\sqcup U_s\sqcup N_s$ such that $|N_s|\leq
(d^{1/(k-1)}-1)(|U_1|+\dots+|U_s|)$. So,
$$
  d^{1/(k-1)}|N_s|\leq (d^{1/(k-1)}-1)(|U_1|+\dots+|U_s|)+(d^{1/(k-1)}-1)|N_s|=|V|(d^{1/(k-1)}-1),
$$
or $\ds |N_s|\leq |V|\frac{d^{1/(k-1)}-1}{d^{1/(k-1)}}$.

Assume now that $|N_s|\leq f(n-2,k)$; then one may color the vertices of~$G(N_s)$ in $n-2$ colors, and then color the vertices of each bipartite graph $G(U_i)$ in two remaining colors. This coloring might be not proper only if some vertices of two subgraphs~$G(U_i)$ and~$G(U_j)$ ($i\neq j$) are adjacent, which is impossible by the condition~(i). So, $G$ is $n$-colorable which is wrong. Therefore, $|N_s|\geq f(n-2,k)+1$ and hence $\ds |V|\geq \frac{d^{1/(k-1)}}{d^{1/(k-1)}-1}|N_s|\geq
\frac{d^{1/(k-1)}}{d^{1/(k-1)}-1}{(f(n-2,k)+1)} $, as desired.
\qed

\Zam 1. In the statement of the Lemma above, one may use the number
$\ds d'=\max_{\emp\ne V'\subseteq V}\min_{u\in V'}|U_{k-1}(u,G(V'))|$ instead of~$d$. For reaching that, on each step it is sufficient to choose the vertex $v\in V_{s-1}$ such that
$$
  |U_{k-1}(v,G_{s-1})|=\min_{u\in V_{s-1}}|U_{k-1}(u,G_{s-1})|.
$$
Clearly, we have $d'\leq d$.

\Zam2.
 On the other hand, the number $d^{1/(k-1)}$ in the same statement can be replaced by $(f(n,k)+1)^{1/k}$. Now, in the proof one may deal with $k+1$ sets
 $$
   U_0(v,G_{s-1})=\{v\},\quad  U_1(v,G_{s-1}), \quad \dots, \quad U_k(v,G_{s-1})
 $$
 and use the condition $|U_k(v,G_{s-1})|\leq |V|=f(n,k)+1$.

\smallskip
The next theorem follows immediately from the Lemmas~\ref{Sergey} and~\ref{Ilya}.

\begin{Theorem}
  \label{recurr}
  For all integer $n,k\geq 2$, we have
  \begin{equation}
    f(n,k)\geq
    \min_{t\geq n(k-1)+1}\max\left\{f(n-1,k)+t,\frac{t^{1/(k-1)}}{t^{1/(k-1)}-1}(f(n-2,k)+1)-1\right\}.
    \label{minmax}
  \end{equation}
  \label{theo}
\end{Theorem}

\proof
From the choice of~$G$ we have $f(n,k)=|G|-1$. From Lemmas~\ref{Sergey} and~\ref{Ilya} it follows that
$$
  |G|\geq \max\left\{f(n-1,k)+d,\frac{d^{1/(k-1)}}{d^{1/(k-1)}-1}(f(n-2,k)+1)-1\right\}+1.
$$
Since $d\geq n(k-1)+1$ by Lemma~\ref{big okr}, the statement holds.
\qed

\begin{Cor}
  For every real $g>1$, we have
  \begin{equation}
    f(n,k)\geq
    \min\left\{f(n-1,k)+g,
    \frac{g^{1/(k-1)}}{g^{1/(k-1)}-1}(f(n-2,k)+1)-1\right\}.
    \label{min}
  \end{equation}
  \label{cor}
\end{Cor}

\proof Let $t_0$ be the integer for which the minimum in~\eqref{minmax} is achieved. As $t>1$ increases,  the value of $f(n-1,k)+t$ also increases, while the value of $\ds \frac{t^{1/(k-1)}}{t^{1/(k-1)}-1}(f(n-2,k)+1)-1$ decreases. Thus, if  $g\leq t_0$, then we have
$$
  {f(n-1,k)+g}\leq {f(n-1,k)+t_0}\leq f(n,k).
$$
Otherwise, we have $g>t_0$ and
$$
  \frac{g^{1/(k-1)}}{g^{1/(k-1)}-1}{(f(n-2,k)+1)}-1\leq
  \frac{t_0^{1/(k-1)}}{t_0^{1/(k-1)}-1}{(f(n-2,k)+1)}-1\leq f(n,k).
  \eqno\qed
$$

\section{Explicit bounds}

Now we present the explicit lower bounds for $f(n,k)$ following from the results of the previous section.

Notice that for every $k$ we have $f(1,k)=1$ and $f(2,k)=2k$. Lemma~\ref{Sergey} implies now the following statement.

\begin{Theorem}
  For all integer $n\geq 1$ and $k\geq 2$ the inequality
  $\ds f(n,k)\geq n+\frac{(k-1)(n-1)(n+2)}2$ holds.
  \label{Sergey2}
\end{Theorem}

\proof Induction on~$n$. In the base cases $n=1$ or $n=2$ the statement holds. Assume now that $n>2$. By Lemmas~\ref{big okr} and~\ref{Sergey} we have $f(n,k)\geq f(n-1,k)+n(k-1)+1$. Next, the hypothesis of the induction implies that
$$
  f(n-1,k)\geq (n-1)+\frac{(k-1)(n-2)(n+1)}2.
$$
Therefore,
$$
f(n,k)\geq f(n-1,k)+n(k-1)+1\geq %n-1+\frac{(k-1)(n-2)(n+1)}2+n(k-1)+1=
n+\frac{(k-1)(n-1)(n+2)}2,
$$
as desired.
\qed

The next estimate uses the whole statement of the Theorem~\ref{theo}. For the convenience, we use the notation $n\^k=n(n+1)\dots(n+k-1)$.

\begin{Lemma}
  Suppose that for some integer $n_0\geq 1$, integer $k\geq 2$, and real $a$, the inequality
  \begin{equation}
    f(m,k)\geq \frac{(m+a)\^k}{2^{k-1}k^k}
    \label{asym1}
  \end{equation}
  holds for two values $m=n_0$ and $m=n_0+1$. Then the same estimate holds for all integer $m\geq n_0$.
  \label{common}
\end{Lemma}

\proof
We prove by induction on $n\geq n_0$ that the estimate~\eqref{asym1} holds for $m=n$. The base cases $m=n_0$ and $m=n_0+1$ follow from the theorem assumptions.

For the induction step, suppose that $n\geq n_0+2$. Let $c=2^{1-k}k^{-k}$, $g=ck(n+a)\^{k-1}$. By the induction hypothesis, we have
\begin{equation}
  f(n-1,k)+g\geq c(n+a-1)\^k+ck(n+a)\^{k-1}
  =c(n+a)\^{k-1}(n+a-1+k)=c(n+a)\^k.
  \label{est-lm}
\end{equation}
Notice that Lemmas~\ref{big okr} and~\ref{Sergey} imply that $f(n,k)\geq
f(n-1,k)+n(k-1)+1$. Hence, if $g\leq n(k-1)+1$, then $f(n,k)\geq
f(n-1,k)+g\geq c(n+a)\^k$, as desired.

Thus we may deal only with the case $g>n(k-1)+1$; in particular, $g>1$. We intend to use Corollary~\ref{cor}; for this, let us estimate the second term in the right-hand part of~\eqref{min}.

From the AM--GM inequality we have
$$
  g^{1/(k-1)}=(ck)^{1/(k-1)}\left((n+a)\^{k-1}\right)^{1/(k-1)}
  \leq \frac1{2k}\left(n+a+\frac k2-1\right).
$$
Let $s=n+a+\frac k2-1$; then $s\geq 2kg^{1/(k-1)}>2k$. Therefore,
\begin{multline*}
  \frac{g^{1/(k-1)}}{g^{1/(k-1)}-1}
  \geq \frac{s}{s-2k}\geq \frac{s+k-1}{s-(k+1)}\geq \\
  \geq \frac{s^2+s(k-1)+\frac{k(k-2)}4}{s^2-s(k+1)+\frac{k(k+2)}4}
  =\frac{(n+a+k-2)(n+a+k-1)}{(n+a-2)(n+a-1)}.
\end{multline*}

Finally, from the induction hypothesis we get
\begin{multline}
  \frac{g^{1/(k-1)}}{g^{1/(k-1)}-1}(f(n-2,k)+1)-1
  \geq \frac{g^{1/(k-1)}}{g^{1/(k-1)}-1}f(n-2,k)\geq\\
  \geq \frac{(n+a+k-2)(n+a+k-1)}{(n+a-2)(n+a-1)}\cdot
  c(n+a-2)\^k=c(n+a)\^k.
  \label{est-lm2}
\end{multline}
Thus, for the value of~$g$ chosen above, Corollary~\ref{cor} and the estimates~\eqref{est-lm} and~\eqref{est-lm2} provide that
$$
  f(n,k)\geq
  \min\left\{f(n-1,k)+g,
    \frac{g^{1/(k-1)}}{g^{1/(k-1)}-1}(f(n-2,k)+1)-1\right\}
  \geq c(n+a)\^k,
$$
as desired.
\qed

Finally, let us show that the constant $a$ in the previous Lemma can be chosen relatively large.

\begin{Theorem}
  For all $k\geq 2$ and $n\geq 2$, we have $\ds f(n,k)\geq \frac{\left(n+k\right)\^k}{2^{k-1}k^k}$.
  \label{estimate}
\end{Theorem}

\proof
Set $a=k$. Let us check the inequality~\eqref{asym1} for $n=2$ and $n=3$. Recall that $f(2,k)=2k$. Now for $m=2$ we get
$$
  2^{k-1}k^kf(2,k)=2^kk^{k+1}=(2k)^{k-1}\cdot 2k^2\geq (k+2)(k+3) \dots
  2k\cdot (2k+1)=(k+2)\^k.
$$
For $m=3$, Theorem~\ref{Sergey2} yields $f(3,k)\geq 5k-2$, and the previous estimate now implies that
$$
  2^{k-1}k^kf(3,k)\geq 2^{k-1}k^k(5k-2)\geq 2^kk^kf(2,k)
  \geq 2(k+2)\^k>(k+3)\^k.
$$
Thus, the inequality~\eqref{asym1} holds for $m=2$ and $n=3$, and hence for all $n\geq 2$ by Lemma~\ref{common}.
\qed

The authors are very grateful to the referees for the valuable comments.

\end{document}